\documentclass{article}

\usepackage {graphicx}
\usepackage {amsfonts}
\usepackage {amsmath}
\usepackage {amssymb}
\usepackage{amsthm}
\theoremstyle{plain}

\usepackage[cp1250]{inputenc}

\newtheorem{proposition}{Proposition}
\newtheorem{remark}{Remark}
\newtheorem{definition}{Definition}
\newtheorem{theorem}{Theorem}
\newtheorem{corollary}{Corollary}
\newtheorem*{mainlemma*}{Main Lemma}
\newtheorem*{lemma*}{Lemma}

\newcommand{\st}[1]{{\cal #1}}
\newcommand{\newton}[2]{{#1 \choose #2}}
\def\blacksquare{\hfill\hbox{\vrule width 5pt height 5pt depth 0pt}}

\begin{document}

\title{Generalized Gaussian processes and relations with random matrices and positive definite functions on permutation groups}

\author{Marek Bo\.{z}ejko$^1$ and Wojciech Bo\.{z}ejko$^2$\\[5mm]
$^1$Institute of Mathematics,\\ University of Wroc\l aw,\\
Pl. Grunwaldzki 2/4,\\ 50-384 Wroc\l aw, Poland\\
email: {\tt marek.bozejko@math.uni.wroc.pl}\\[5mm]
$^2$Institute of Computer Engineering, Control and Robotics\\Wroc\l aw University of Technology,\\
Wyb. Wyspia\'{n}skiego 27,\\ 50-370 Wroc\l aw, Poland\\
email: {\tt wojciech.bozejko@pwr.wroc.pl}}

\date{}

\maketitle
\begin{abstract}
The main purpose of this paper of the paper is an explicite construction of generalized Gaussian process with function $t_b(V)=b^{H(V)}$, where $H(V)=n-h(V)$, $h(V)$ is the number of singletons in a pair-partition $V \in \st{P}_2(2n)$.

This gives another proof of Theorem of A. Buchholtz [Buch] that $t_b$ is positive definite function on the set of all pair-partitions.

Some new combinatorial formulas are also presented. Connections with free additive convolutions probability measure on $\mathbb{R}$ are also done. Also new positive definite functions on permutations are presented and also it is proved that the function $H$ is norm (on the group $S(\infty)=\bigcup S(n)$.
\end{abstract}

\noindent {\it MSC:} primary 46L54; secondary 05C30\\

\noindent {\bf Keywords:} Infinite divisibility; Gaussian processes; Pair-partitions; Positive definite functions; Permutation group; Random matrices; Free convolutions

\section{Introduction}

We present some new construction of generalized Gaussian processes and its relations with random matrices as well as with positive definite functions defined on permutations groups. The plan of the paper is following: first we present definitions and remarks on pair-partitions. Next, Markov random matrices and function $h$ on pair-partitions are presented in the Section 3 -- $\st{P}_2(2n)$ as obtained by Bryc, Dembo, Jiang \cite{B-D-J}.
Generalized strong Gaussian processes (fields) $\{G(f)$ are showed in the Section4, $f \in \st{H}\}$ (GSGP), $\st{H}$ - real Hilbert space, as well as the main and the new examples. The main theorem is placed in the Section~5. The free product of (GSGP) $G(f)$ and free Gaussian field $G_0(f)$ is again (GSGP) -- this appears in this paper as the Theorem~3. In the Section~5 we also present the new interesting results on free convolutions of measures extending results of Bo\.{z}ejko and Speicher \cite{B-Sp2}. Such negative-definite functions have mane applications in the field of telecommunication, parallel and quantum computing as well as in operations research (see \cite{B-Wo,Boz2}). We show that the function $H(\sigma) = n- h(\sigma)$, $\sigma \in S(n)$ is conditionally negative definite, i.e. for each $x>0$, $\exp(-x\;H(\sigma))$ is positive functions on the permutation group $S(\infty)=\bigcup S(n)$. Also in the Theorem~6 (Section 6) it is shown that the function $H$ defined as $d_H(\sigma,\tau)=H(\sigma^{-1}\tau)$ is left-invariant distance on~$S(\infty)$.

\section{Definitions and remarks on pair-partitions}

\begin{definition}
Let $\gamma_0$ is the Wigner (semicircle) law with density $\frac{1}{\sqrt{2 \pi}} \sqrt{4-x^2} dx$. This is free Gaussian law. By
$\gamma_1$ we denote the Normal law $N(0,1)$ with density $\frac{1}{\sqrt{2 \pi}} e^{-x^2/2}$. The moments (only even)
\begin{equation}
\int_{-\infty}^\infty x^{2n} d \gamma_1 (x) = m_{2n}(\gamma_1) = \sum_{V \in \st{P}_2 (2n)} 1 = 1 \cdot 3 \cdot 5 \cdot \ldots \cdot (2n-1) =
(2n-1)!!=p_{2n},
\end{equation}
where $\st{P}_2(2n)$ is the set of all pair-partitions on $2n$-elementary set $\{1, 2, \ldots, 2n\}$.
\end{definition}
The moments of the Wigner law
$m_{2n}(\gamma_0)=\frac{1}{n+1} \newton{2n}{n}=$ the cardinality of all non-crossing pair-partitions of $\st{P}_2(2n)$,
where a partition $V \in \st{P}_2(2n)$ has a \emph{crossing} if there exists blocks $(i_1,j_1),(i_2,j_2) \in V$ such that $i_1 < i_2 < j_1 < j_2$:
\begin{center}\begin{minipage}{15mm}\includegraphics[width=15mm]{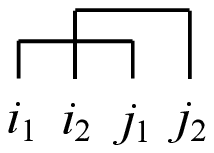}.\end{minipage}\end{center}

\noindent Let $cr(V)=\#$ of all crossings of pair-partition $V$.

In contrary, the partition $V$ is called \emph{non-crossing}; $NC_2(2n)$ -- denote the set of all non-crossing pair-partitions on $2n$-elements set
$\{1,2,\ldots,2n\}$.
Pictorially:

\begin{center}
\begin{minipage}{3cm}\includegraphics[width=3cm]{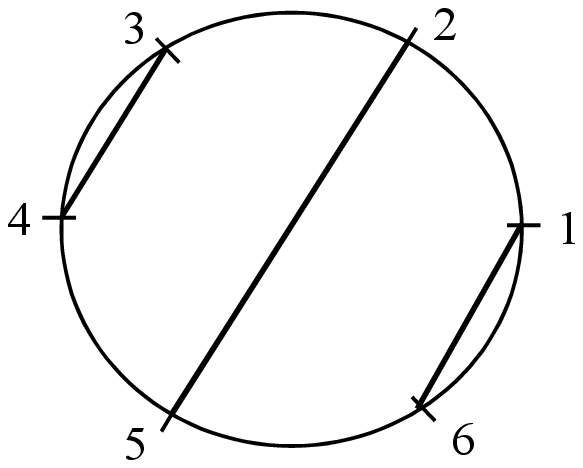}\end{minipage}\hspace{5mm}\begin{minipage}{65mm} $V=\{(1,6),(2,5),(3,4)\}$ or
\includegraphics[width=2cm]{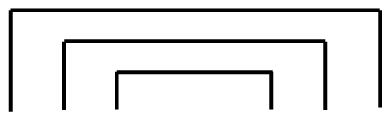} .\end{minipage}
\end{center}

\noindent It is well known, that the cardinality of $NC_2(2n)=\frac{1}{n+1} \newton{2n}{n}$.

\begin{definition}
The block $B \in V \in \st{P}_2(2n)$ is a {\bf singleton}, if $B$ has no crossing with other block $C \in V$.
\end{definition}
Ex. \begin{minipage}{20mm}\includegraphics[width=20mm]{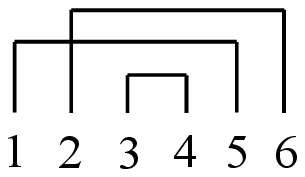}\end{minipage}, the block $(3,4)=B$ is a singleton.\\

\noindent Following \cite{B-D-J}, let us denote $h(V)=\#$ of singletons in the pair-partition $V \in \st{P}_2(2n)$.

\noindent For example: if $V=$ \begin{minipage}{17mm}\includegraphics[width=17mm]{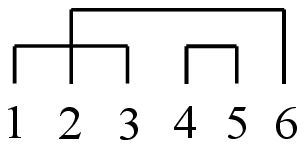}\end{minipage} then $h(V)=1$.

\noindent Facts:
\begin{enumerate}
\item[($\alpha$)] If $V$ is {\bf non-crossing pair-partition} on $\{1,2\ldots,2n\}$, then $h(V)=n.$
\end{enumerate}
The important fact for us is the following:
\begin{enumerate}
\item[($\beta$)] If $V$ is connected and $V \in \st{P}_2(2n), n>1$, then $h(V)=0$,\\[2mm] for example $h($\begin{minipage}{16mm}\includegraphics[width=17mm]{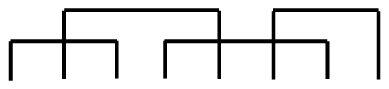}\end{minipage}$)=0$.\\
\end{enumerate}

\noindent The pair-partition $V \in \st{P}_2(2n)$ is \emph{connected} if its graph is \emph{connected} set. For example \begin{minipage}{5mm}\includegraphics[width=5mm]{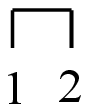}\end{minipage} is connected.
Also \begin{minipage}{17mm}\includegraphics[width=17mm]{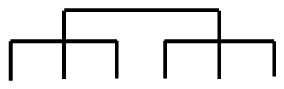}\end{minipage} is connected.\\

\noindent Let $cc(V) =$ the number of connected components of the graph of the partition $V \in \st{P}_2(2n)$, and

\noindent let $c_{2n} = \# \{V \in \st{P}_2(2n): V \mbox{ is connected}\}$.

\noindent That sequence is \emph{the free cumulant} of the classical Gaussian distribution $N(0,1)$, i.e.:
$$c_2=1,\quad c_4=1,\quad c_6=4,\quad c_8=27,\quad c_{10}=248,\ldots.$$
The following formula is due to Riordan \cite{R}, see also Belinschi, Bo\.{z}ejko, Lehner, Speicher \cite{B-B-L-S}:
\begin{equation}\label{eq1}
c_{2n+2}=n \sum_{i=1}^{n} c_{2i} c_{2(n-i-1)}
\end{equation}
and that sequence is the moment sequence of some symmetric probability measure on real line, as it was proved in \cite{B-B-L-S}.

That fact is equivalent to the following result:

\begin{theorem}{\cite{B-B-L-S}}
Normal law $\gamma_1$ is infinitely divisible in free $\boxplus$-convolution (i.e.: $\gamma_1 \in ID(\boxplus)$).
\end{theorem}

One of our aim of this work is to find different proof of that above  result using different method and the function
$h(V)$.

\noindent Let us first prove the following:

\begin{proposition}\label{prop_1}
Let us define $T_{2n}=\sum_{V \in \st{P}_2(2n)} h(V)$, and
$$
p_{2n}=(2n-1)!!=\sum_{V \in \st{P}_2(2n)} 1,\quad p_{2n}=1,3,15,105,\ldots,
$$
then
\begin{equation}\label{eq2}
T_{2n+2}=(n+1) \sum_{k=0}^n p_{2k} \cdot p_{(2n-2k)}.
\end{equation}
\end{proposition}
That sequence $T_{2n}$ is following: 1, 4, 21, 144, 1245, 13140, 164745,\ldots.

\noindent The proof of the formula (\ref{eq2}) is by a simple considerations, if we consider pair-partitions as lying on a circle.

\noindent For example ($n=3$):
\begin{center}
\begin{minipage}{35mm}\includegraphics[width=35mm]{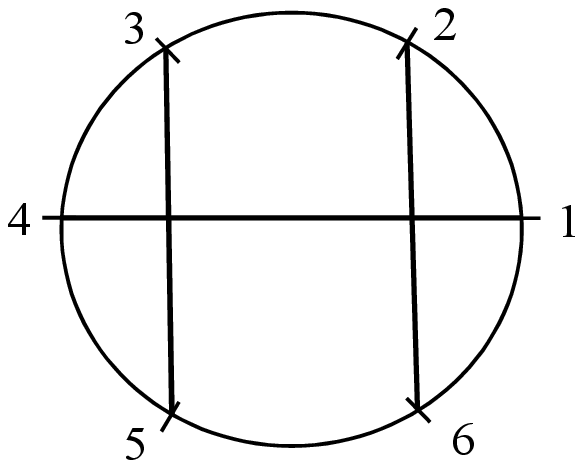}\end{minipage}\hspace{5mm}$V=\{(1,4),(2,6),(3,5)\}$,
\end{center}
\begin{center}
\begin{minipage}{35mm}\includegraphics[width=35mm]{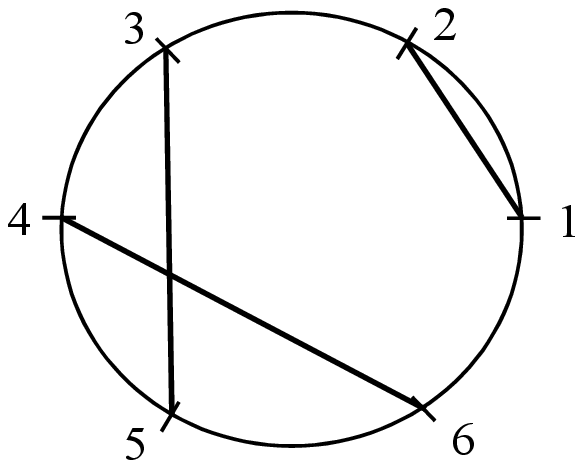}\end{minipage}\hspace{5mm}$V=\{(1,2),(3,5),(4,6)\}$.
\end{center}

\noindent {\bf Proof of Proposition~\ref{prop_1}.} To obtain the proof of the formula (\ref{eq2}), we consider the set of 2-pairing of the set $\{1,2,3,\ldots,2n+2\}$.

Let $A^{(1)}_k=$ all pair-partitions which contains the singletons starting from 1 to $2k$: $(1,2k)$, $k=1,2,\ldots,$ $2n+2$.

$A^{(3)}_k$ the set of $\st{P}_2(2n+2)$, which contains the singletons $(3,4), (3,6), \ldots,$ i.e. this the rotation of the set $A^{(1)}_k$; i.e.

$A^{(1)}_1$\hspace{5mm}\begin{minipage}{35mm}\includegraphics[width=35mm]{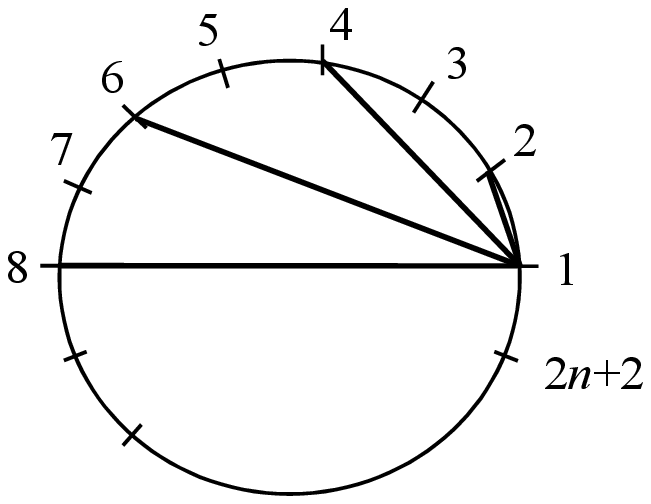}\end{minipage}

Similarly, let $A^{(2j+1)}_k$ is the configuration of 2-pairing like $A^{(1)}$, but starting from the point $2j+1$, $j=0,1,2,\ldots,n$. All the sets $A^{(2j+1)}_k$ are \emph{disjoint}, and the number of singletons in that set  $A^{(2j+1)}_k$ equals $\sum_{k=0}^n p_{2k} \cdot p_{(2n-2k)}$.

Therefore the number of all singletons in $\st{P}_2 (2n+2)$ is equal
$$
T_{2n+2}=(n+1) \sum_{k=0}^n p_{2k} \cdot p_{(2n-2k)}. \hspace{5mm}\blacksquare
$$

\section{Markov random matrices and function $h$ on pair-partitions $\st{P}_2(2n)$}

Let $\{X_{ij} : j \geq i \geq 1 \}$ be an infinite upper triangular array of i.i.d. random variables and define $X_{ji} = X_{ij}$ for $j > i \geq 1$.

\noindent Let $X_n = [ X_{ij} ]_{1 \leq i,j \leq n}$, and
$$
D_n = \mbox{diag}(\sum_{j=1}^n X_{ij})
$$
is a diagonal matrix . We define \emph{Markov matrices} $M_n$ as a random matrix given by
$$
M_n = X_n - D_n,
$$
so then each of rows (and columns) of $M_n$ has a zero sum.

Here for a symmetric $n \times n$ matrix $A$, its empirical distribution is done as
$$
\hat{\mu}(A) = \frac{1}{n} \sum_{j=1}^n \delta_{\lambda_j(A)},
$$
$\lambda_j(A)$, $1 \leq i \leq n$, denote the eigenvalues of the matrix $A$ and $\delta_s$ is the Dirac mass at the point $s \in \mathbb{R}$.

The Theorem of Bryc, Dembo, Jiang \cite{B-D-J} is following:
\begin{theorem}[\cite{B-D-J}]
If $X_{ij}$ are i.i.d. random variables with $\mathbb{E} X_{ij}=0$ and $\mathbb{E} X_{ij}^2=1$, then
$\hat{\mu}(\frac{M_n}{\sqrt{n}})$ converges weakly as $n \rightarrow \infty$ to the measure $\gamma_M= \gamma_0 \boxplus \gamma_1$, where $\boxplus$ denotes the free additive convolution of probability measures.
\end{theorem}

The even moments of the measure $\gamma_M= \gamma_0 \boxplus \gamma_1$ are following:
$$
m_{2n}(\gamma_M) = \sum_{V \in \st{P}_2(2n)} 2^{h(V)}, \mbox{  and  } m_{2n+1}(\gamma_M)=0,
$$
$\gamma_0$ is the Wigner (semicircle) law done by density $\frac{1}{2 \pi} \sqrt{4-x^2} \cdot \chi [-2,2]$, $\gamma_1$ is the Normal law $N(0,1)$.

\section{Generalized Gaussian process (field) $G(f)$, $f \in \st{H}$, $\st{H}$ - real Hilbert space. Main and new examples}

Let in some probability system $(\st{A},\varepsilon)$, (\st{A} - $\ast$-algebra with unit, and $\varepsilon$ is state on $\st{A}$). The
family $G(f)=G(f)^* \in \st{A}$, $f \in \st{H}$ - real Hilbert space, is called \emph{normalized generalized Gaussian process} (GGP), if for each orthogonal map
$\st{O}:\st{H} \rightarrow \st{H}$, for $f_j \in \st{H}_\mathbb{R}$, we have
$$
\varepsilon(G(f_1)G(f_2) \ldots G(f_k))=\varepsilon(G(\st{O}(f_1))G(\st{O}(f_2)) \ldots G(\st{O}(f_k)))=
$$
$$
=\left\{
\begin{array}{ll}
0 & \mbox{for k - odd},\\
\displaystyle\sum_{V \in \st{P}_2(2n)} t(V) \prod_{(i,j) \in V} < f_i | f_j>& \mbox{for $k=2n$},\\
\end{array}
\right.
$$
for some function $t:\st{P}_2(\infty) \rightarrow \mathbb{C}$ which will be called positive definite on $\st{P}_2(\infty)$(see \cite{B-Sp2}, Gu\c{t}\v{a}-Maassen \cite{G-M1,G-M2}, for further facts), with normalization $t(\includegraphics[width=6mm]{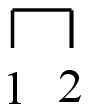})=1$.
The main examples of (GGP) are related to q-CCR relations \cite{B-Sp1},\cite{B-K-S}:\\[4mm]
$\mbox{(q-CCR)} \quad\quad\quad\quad\quad\quad\quad a(f)a^*(q)-q\;a^*(q)a(f)=<f,q>{\rm I}$\\[4mm]
$-1 \leq q \leq 1$, $a(f)\Omega=0$, $f,q \in \st{H}_\mathbb{R}$ and $\| \Omega \| = 1$, $\Omega$ - vacuum vector.

\noindent If we take $G_q(f)=a(f)+a^*(f)$, and as a state -- vacuum state:\\ $\varepsilon(T)=<T \Omega | \Omega>$, then
we get $q$-Gaussian field and $t_q(V)=q^{cr(V)}$ for $V \in \st{P}_2(2n)$, where $cr(V)$ is the number of crossings in a partition $V$.

Other examples were constructed by Bo\.{z}ejko-Speicher \cite{B-Sp2} by the function:
$$
t_s(V)=s^{n-cc(V)}, \quad a \leq s \leq 1.
$$
That examples were important to prove that Normal law $\gamma_1$ is free infinitely divisible, i.e. $\gamma \in ID(\boxplus)$.

\noindent Namely, the following fact was proven in \cite{B-Sp2}. For
$s \geq 1$
$$
m_{2n}(\gamma_1^{\boxplus s}) = \sum_{V \in \st{P}_2(2n)} s^{cc(V)},
$$
where $\boxplus$ is the free additive convolution.
Many other examples were done by Accardi-Bo\.{z}ejko \cite{A-B}, Gu\c{t}\v{a}-Maassen \cite{G-M1,G-M2}, Bo\.{z}ejko-Yoshida \cite{B-Y}, Bo\.{z}ejko-Gu\c{t}\v{a} \cite{B-G}, Bo\.{z}ejko \cite{Boz1} and Bo\.{z}ejko-Wysocza\'{n}ski \cite{B-W}.

Our work presents among others the proof of the result of A. Buchholz \cite{Buch}, that it exists a explicite realization of generalized Gaussian process connected with the function of Bryc-Dembo-Jiang
$$
t_b(V)=b^{n-h(V)}=b^{H(V)}
$$
for $0 \leq b \leq 1$, $V \in \st{P}_2(2n)$, and this is consequence of \emph{our Main Theorem}.

We define generalized \underline{strongly} Gaussian process $(G(f),t_G,\varepsilon)$, $f \in \st{H}$, as generalized \emph{Gaussian process}, such that the function $t_G = t$ on $\st{P}_2(\infty)$ is
\emph{strongly multiplicative}, i.e.
$$
t(V_1 \dot{\cup} V_2) = t(V_1) \cdot t(V_2),
$$
for $V_1, V_2 \in \st{P}_2(\infty)$, which are pair-partitions on disjoint sets.

\noindent The simple examples of strongly \emph{multiplicative} functions on $\st{P}_2(\infty)$ are following:\\
$$
\begin{array}{ll}
1^{\rm o} & q^{cr(V)},\\
2^{\rm o} & s^{n-cc(V)},\\
3^{\rm o} & b^{n-h(V)},\\
\end{array}
$$
see also \cite{B-G} for more strongly multiplicative examples, related to the Thoma characters on $S(\infty)$ group.

\noindent That classes of processes correspond to so called \emph{pyramidal} independence, which has been considered by B. K\"{u}mmerer \cite{K}, see also \cite{B-Sp2}.\\

\noindent In all above examples: if $q=b=s=1$, we get classical Gaussian process.\\

\noindent Important fact: if $\st{H} = L^2 (\mathbb{R}^+,dx)$ and $B_u = G(\chi_{[0,u]})$, where $\chi_{[0,u]}$ is the characteristic function of interval $[0,u]$, then $B_u$, for $u \geq 0$, is a realization of classical \emph{Brownian motion}.\\

\noindent If we take $q=b=s=0$, we get a construction of the \emph{free Brownian motion} (see books of Voiculescu, Dykema, Nica \cite{V-D-N}, Nica, Speicher, \cite{N-S} and Hiai-Petz \cite{H-P}).\\

\section{The main theorem}

\begin{theorem}\label{th3}
If $G(f)$ is normalized generalized strong Gaussian process, $f \in \st{H}$, $\st{H}$ is a real Hilbert space, and $G_0(f)$ is the free Gaussian process, and
operators $\{G(f): f \in \st{H}\}$ and $\{G_0(f): f \in \st{H}\}$ are \emph{free independent} in some probability system $(\st{A},\varepsilon)$, then
for each $0 \leq b \leq 1$, the process:
$$
X_b(f) = \sqrt{b}\; G(f) + \sqrt{1-b}\; G_0(f),\quad f \in \st{H}
$$
is again generalized \underline{strong} Gaussian process.\\

Moreover, if $G(f)$ corresponds to \emph{strongly multiplicative} function
$$
t_G: \bigcup_{\st{A}_n} \st{P}_2(2n)\rightarrow\mathbb{C}
$$
done by equation
$$
\varepsilon(G(f_1) \ldots G(f_k)) = \left\{
\begin{array}{ll}
\sum_{V \in \st{P}_2(2n)} t_G(V) \prod_{(i,j) \in V} <f_i | f_j >, & \mbox{if $k=2n$},\\
0, & \mbox{if $k$ odd}.\\
\end{array}
\right.
$$

Then the corresponding \emph{strongly multiplicative} function of the generalized Gaussian process $X_b$, $0 \leq b \leq 1$ is following:
$$
t_{X_b}(V)=b^{H(V)} \cdot t_G(V),\quad \mbox{for $V \in \st{P}_2(2n)$.}
$$
i.e.
$$
\varepsilon(X_b(f_1)\; X_b(f_2)\ldots X_b(f_{2n})) = \sum_{V \in \st{P}_2(2n)} b^{H(V)} t_{G(V)} \prod_{(i,j) \in V} <f_i | f_j >,
$$
for $f_i \in \st{H}$, and odd moments are zero.
\end{theorem}

\noindent In the proof of the Theorem~\ref{th3} we will need the following Lemma.

\begin{mainlemma*} Let t be \emph{strongly multiplicative} function corresponding to strongly generalized Gaussian process (field) $G(f)=G_t(f)$, $f \in \st{H}_{\mathbb{R}}$ such that
$$
\varepsilon(G(f_1) \ldots G(f_k)) = \left\{
\begin{array}{ll}
\displaystyle\sum_{V \in \st{P}_2(2n)} t(V) \displaystyle\prod_{(i,j) \in V} <f_i | f_j >, & \mbox{$k=2n$},\\
0, & \mbox{if $k=2n+1$}.\\
\end{array}
\right.
$$
then the free cumulants are following:
\begin{equation}\label{eq_asterisk}
r_{k}(G(f_1),G(f_2),\ldots,G(f_k)) =
\left\{
\begin{array}{ll}
\hspace*{-2mm}\displaystyle\sum_{\scriptsize
\begin{array}{c}
V \in \st{P}_2(2n)\\
cc(V)=1\\
\end{array}
} \hspace*{-5mm}t(V) \displaystyle\prod_{(i,j) \in V} <f_i | f_j>, & k=2n,\\
0, & k=2n+1.
\end{array}
\right.
\end{equation}
\end{mainlemma*}

\noindent {\bf Proof.} Let $NC_e(2n)$ denotes the set of all even non-crossing partition $\st{V}$ of $2n$, which all blocks of $\st{V}$ are even.

As in the proof of Theorem~11 in \cite{B-Y}, we define the mapping $\Phi:\st{P}_2(2n) \rightarrow NC_e(2n)$ as follows: given a pair-partition $V \in \st{P}_2(2n)$,
the connected components of $V$, will induce the even non-crossing partition $\Phi(V)=W$.

For example, if
\begin{center}
\includegraphics[width=88mm]{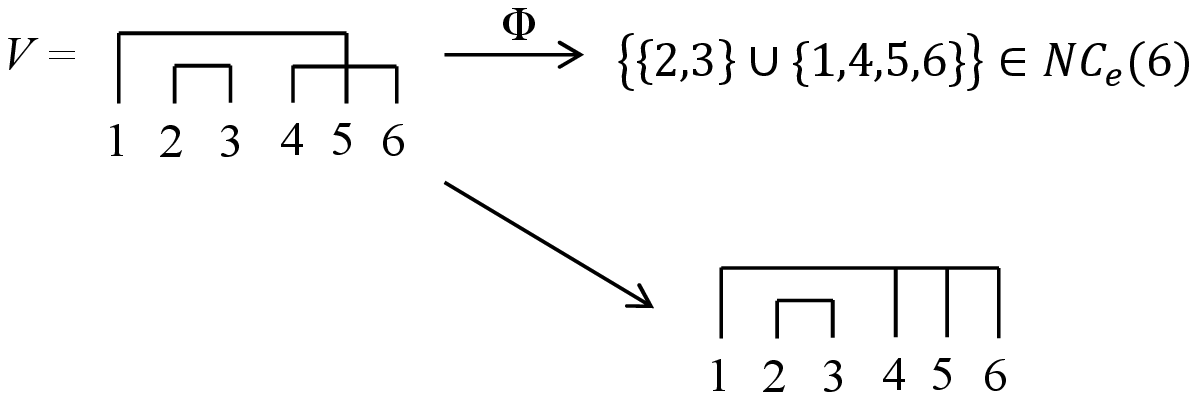}
\end{center}
(i.e. $\Phi$ in some sense forgets crossings of partitions).

Let us denote $G(f_j)=g_j$ in the proof. Since all odd moments of $g_j$ vanish, hence also all odd free cumulants are vanish, i.e.
$r_{2k+1}(g_{i_1},\ldots,g_{i_{2k+1}})=0$.

Therefore the free moment-cumulant formula is following:
\begin{equation}\label{eq_gwiazdka}
\varepsilon(g_1 g_2 \ldots g_{2n}) =
\displaystyle\sum_{V \in NC_e(2k)} \displaystyle\prod_{\scriptsize
\begin{array}{c}
B \in V\\
B=\{i_1,\ldots,i_{2s}\}\\
\end{array}
} r_{2s} (g_{i_1},\ldots g_{i_{2s}}).
\end{equation}

By the assumption we have
\begin{equation}\label{eq_2gwiazdki}
\varepsilon(G(f_1)\ldots G(f_{2n})) = \varepsilon(g_1 g_2 \ldots g_{2n})=
\displaystyle\sum_{V \in \st{P}_2(2n)} t(V) \displaystyle\prod_{(i,j) \in V} <f_i|f_j>
\end{equation}
where $t$ is strongly multiplicative.

Let us denote
$$
\tilde{r}_{2k}(g_1,g_2,\ldots,g_{2k})=
\displaystyle\sum_{\scriptsize
\begin{array}{c}
\rho \in \st{P}_2(2k)\\
cc(\rho)=1\\
\end{array}
} t(\rho) \prod_{(i,j) \in \rho} <f_i|f_j>.
$$
We want to show that
$$
\tilde{r}_{2k} (g_1,\ldots,g_{2k})=r_{2k} (g_1,\ldots,g_{2k}).
$$
By the strong multiplicativity of $t$, we have that the function
$$\tilde{t}(\pi)=t(\pi)\prod_{(i,j)\in \pi} <f_i|f_j>$$
is also strong multiplicative on the $\st{P}_2(2n)$.

Let us fix a non-crossing partition $V \in NC_e(2n)$. By the strong multiplicative property of $\tilde{t}$, we have
$$
\displaystyle\sum_{\scriptsize
\begin{array}{c}
\pi \in \st{P}_2(2n)\\
\Phi(\pi)=V\\
\end{array}
} \tilde{t}(\pi)=
\hspace*{-7mm}\displaystyle\prod_{\scriptsize
\begin{array}{c}
B \in V\\
B=\{i_1,i_2,\ldots,i_{2s}\}\\
\end{array}
} \hspace*{-7mm}\tilde{r}_{2s} (g_{i_1},g_{i_2},\ldots,g_{i_{2s}}).
$$
Therefore the formulas (\ref{eq_asterisk}) and (\ref{eq_gwiazdka}) implies that
$$
\varepsilon(g_1 g_2 \ldots g_{2n})=
\displaystyle\sum_{V \in NC_e(2n)}
\hspace*{-2mm}\displaystyle\prod_{\scriptsize
\begin{array}{c}
B \in V\\
B=\{i_1,i_2,\ldots,i_{2s}\}\\
\end{array}
} \hspace*{-7mm}\tilde{r}_{2s} (g_{i_1},g_{i_2},\ldots,g_{i_{2s}}).
$$
Comparing with the formulas (\ref{eq_asterisk}) and (\ref{eq_gwiazdka}) and using M\"{o}bius inversion formula for the lattice of non-crossing partition (see Nica, Speicher book \cite{N-S}) we get
$$
\tilde{r}_{2s} (g_1,\ldots,g_{2s}) = r_{2s} (g_1,\ldots,g_{2s}).
$$\hspace{5mm}\blacksquare

\noindent Now we can start the proof of the Main Theorem using the Main Lemma.\\

\noindent {\bf Proof of the Main Theorem.}  By definition of the freeness of the families $\{ G(f) \}_{f \in \st{H}}$, $\{ G_0 (f) \}_{f \in \st{H}}$ in the probability system $(\st{A},\varepsilon)$, we have that all mixed free cumulants
$$
r_k (G_{\varepsilon_1} (f_1), G_{\varepsilon_2} (f_2),\ldots,G_{\varepsilon_k} (f_k))=0,\; \mbox{for all}\; k=2,3,\ldots,
$$
if the sequence $(\varepsilon_1,\varepsilon_2,\ldots,\varepsilon_k)$ is not constant ($\epsilon_j \in \{0,1\})$ and in the proof we denote $G_1 (f) \stackrel{def}= G(f)$.

Therefore the free cumulants of $X_b (f)=\sqrt{b} \; G(f) + \sqrt{1-b} \; G_0 (f)$ are following:
$$
r_{2k} ( X_b (f_1),\ldots,X_b (f_{2k}))=
$$
\begin{equation}\label{eq_proof_main_th}
=b^k \; r_{2k} (G(f_1),\ldots,G(f_{2k}))+(1-b)^k \; r_{2k} (G_0 (f_1),\ldots,G_0 (f_{2k}))
\end{equation}
and all odd free cumulants of $X_b (f)$ are zero.

From the assumption $G_0 (f)$ is the free normalized Gaussian process, i.e.
$$
t_{G_0}(V)=\left\{
\begin{array}{ll}
1, & \mbox{if $V \in NC_2(2k)$},\\
0, & \mbox{otherwise}.
\end{array}
\right.
$$
and
this is equivalent that
$$
r_2 (G_0 (f_1), G_0 (f_2))=<f_1,f_2>
$$
and
$$
r_{2k} (G_0 (f_1),\ldots,G_0 (f_{2k})) = 0,
$$
for $k>1$ and arbitrary $f_j \in \st{H}_{\mathbb{R}}$.

That above facts follow at once from our Main Lemma, as we can see now:

Since by definition:
$$
\varepsilon(G_0 (f_1),\ldots,G_0 (f_{2k}))=\sum_{V \in NC_2(2k)} \prod_{(i,j) \in V} <f_i|f_j>
$$
hence by Main Lemma:
$$
r_{2k}(G_0 (f_1),\ldots,G_0 (f_{2k}))=\sum_{
\scriptsize
\begin{array}{c}
V \in NC_2(2k)\\
cc(V)=1\\
\end{array}
} \prod_{(i,j) \in V} <f_i|f_j>.
$$
But if $V \in NC_2(2k)$, and $cc(V)=1$, then we have that $k=1$ and therefore for $k>1$
$$
r_{2k}(G_0 (f_1),\ldots,G_0 (f_{2k}))=0.
$$

If now $V \in \st{P}_2(2k)$ and $cc(V)=1$, $k>1$, then we have $H(V)=k-h(V)=k$, we get by (\ref{eq_proof_main_th}) and Main Lemma
$$
r_{2k}( X_b (f_1),\ldots,X_b (f_{2k}) ) = b^k \; r_{2k} (G(f_1),\ldots,G(f_{2k}))=
$$
$$
= \hspace{-5mm}\sum_{
\scriptsize
\begin{array}{c}
V \in \st{P}_2(2k)\\
cc(V)=1\\
\end{array}
}
b^k \; t(V)
\prod_{(i,j) \in V} <f_i|f_j>=\hspace{-5mm}
\sum_{
\scriptsize
\begin{array}{c}
V \in \st{P}_2(2k)\\
cc(V)=1\\
\end{array}
}
b^{H(V)} \; t(V)
\prod_{(i,j) \in V} <f_i|f_j>.
$$
On the other hand, for $k=1$ we have
$$
r_2(X_b(f_1),X_b(f_2))=b<f_1|f_2>+(1-b)<f_1|f_2>=<f_1|f_2>,
$$
since $\{G(f)\}$ is normalized Gaussian process, i.e. $t($\includegraphics[width=5mm]{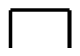}$)=1$
and since\linebreak $H($\includegraphics[width=5mm]{fig18.eps}$)=1-h($\includegraphics[width=5mm]{fig18.eps}$)=0$.

Therefore for all $k \geq 1$
$$
r_{2k}( X_b (f_1),\ldots,X_b (f_{2k}))=
\hspace{-5mm}
\sum_{
\scriptsize
\begin{array}{c}
V \in \st{P}_2(2k)\\
cc(V)=1\\
\end{array}
}
b^{H(V)} \; t(V)
\prod_{(i,j) \in V} <f_i|f_j>
$$

Now using again the free moment-cumulant formula, our Main Lemma and the strong multiplicativity of the function
$b^{H(V)}t(V)$, we get
$$
\varepsilon(X_b(f_1),X_b(f_2),\ldots,X_b(f_{2k}))=\sum_{V \in \st{P}_2(2k)} b^{H(V)}\; t(V) \prod_{(i,j) \in V} <f_i|f_j>.
\hspace{3mm}\blacksquare$$\\

\noindent After that considerations a natural problem appears:\\

\noindent {\bf Problem 1.} Is it true that if we have 2 strongly generalized Gaussian processes $\{G_1(f)\}_{f \in \st{H}_{\mathbb{R}}}$,
$\{G_2(f)\}_{f \in \st{H}_{\mathbb{R}}}$ which are free, that the Gaussian process $Z(f)=G_1(f)+G_2(f)$ is again
strongly Gaussian?\\

As a corollary from the Main Lemma we get well-known similar simple proposition for probability measures on real line (see \cite{B-Sp2}, \cite{B-D-J}, \cite{Leh}):

\begin{proposition}
If $\mu$ is symmetric measure on $\mathbb{R}$ with all moments, and
$$
m_{2n}(\mu)= \sum_{V \in \st{P}_2(2n)} t(V), \mbox{ where $t$ -- strongly multiplicative,}
$$
 then the free cumulants of the measure $\mu$ are of the form:
 $$
 r_{2n}(\mu)=\sum_{\scriptsize
\begin{array}{c}
V \in \st{P}_2(2n)\\
cc(V)=1\\
\end{array}
} \hspace*{-4mm} t(V).
 $$

\end{proposition}

\noindent Now we show some special case of our results.\\

In particular case, let $f \in \st{H}$, $\| f \| = 1$, and let the law of $G(f)$ is the probability measure $\mu$ on $\mathbb{R}$, \quad $\st{L}(G(f))=\mu$, i.e.
$$
\varepsilon(G(f)^k) = \int_\mathbb{R} \lambda^k d \mu(\lambda), \quad k=0,1,2,\ldots
$$
and let $\gamma_0$ is the law of the Wigner-semicircle-free Gaussian law $G_0(f)$, then the law of the process $X_b(f)$:
$$
\st{L}(X_b(f))=D_{\sqrt{b}}(\mu) \boxplus D_{\sqrt{1-b}}(\gamma_0) = \mu_b,
$$
here $\boxplus$ is the \emph{free additive convolution}, and $D_\lambda$ is \emph{the dilation of the measure done by the formula:}
$$
(D_\lambda \mu)(E) = \mu(\lambda^{-1} E),\; \mbox{for Borel set}\; E \subset \mathbb{R},\quad \lambda > 0.
$$
Hence from the Main Theorem we get that the even moments of the measure $\mu_b$ are following:
$$
m_{2n}(\mu_b)=\int \lambda^{2n} d \mu_b(\lambda) = \sum_{V \in \st{P}_2(2n)} b^{H(V)} \cdot t_G (V)
$$
and $m_{2n+1} (\mu_b)=0$.

If we take the classical Gaussian process as $G(f)$, corresponding to $t(V) \equiv 1$, for all $V \in \st{P}_2(2n)$, se we get as corollary the completely different proof of the theorems of A. Buchholz \cite{Buch}.

\begin{corollary}\label{corollary_main_theorem}
For all $0 \leq b \leq 1$, the function $t(V)=b^{H(V)}$ is \emph{strongly multiplicative}, \emph{tracable} and \underline{positive definite} on the set of all pair-partitions $\st{P}_2(\infty)=\bigcup_n \st{P}_2(2n)$.
\end{corollary}

Here the function $H(V)=n-h(V)$ is \emph{tracable}, i.e.
$$
H(V)=H(\overrightarrow{V}),\;\; \mbox{where}
$$
for a pair-partition $V=\{(i_1,j_1),(i_2,j_2),\ldots,(i_n,j_n)\}$,
$\overrightarrow{V}=\{(\overline{i}_1,\overline{j}_1),(\overline{i}_2,\overline{j}_2),\ldots$,
$(\overline{i}_n,\overline{j}_n)\}$, $\overrightarrow{V}$ is the cyclic rotation of our partition $V$
($i_k \rightarrow 1+i_k (\mbox{ modulo } 2n)$).

\noindent For example:\\
\hspace*{31mm}$V$\hspace{53mm}$\overleftarrow{V}$
\begin{center}
\includegraphics[width=90mm]{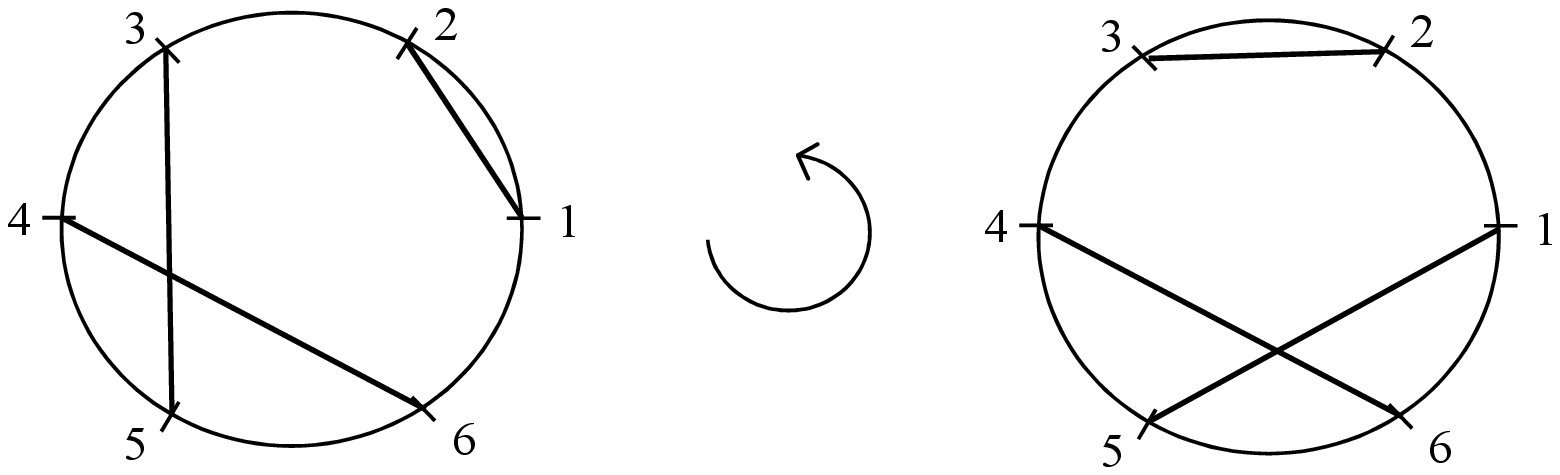}
\end{center}

\begin{remark}
If $(G(f),t,f \in \st{H})$ is generalized Gaussian process and $t$ is {\bf tracable}, then
$\varepsilon$ is a \emph{trace} on the algebra generated by $G(f)$, $f \in \st{H}$, i.e.
$$
\varepsilon(G(f_1) \ldots G(f_k))=\varepsilon(G(f_k) G(f_1) \ldots G(f_{k-1})),
$$
in general $\varepsilon(XY)=\varepsilon(YX)$, for $X,Y$ in $\ast$-algebra generated by the field $G(f)$, $f \in \st{H}$.
\end{remark}

\section{Free convolutions of measures}

As Corollary~\ref{corollary_main_theorem} from the \emph{Main Theorem} we get the following generalization of the Theorem~6 from [B-Sp2].

\begin{proposition}
Let for $0 \leq b \leq 1$, $\rho_b = D_{\sqrt{b}}(\gamma_1) \boxplus D_{\sqrt{1-b}}(\gamma_0)$ then for
$0 \leq c \leq 1$
\begin{equation}\label{eq3}
\rho_{(bc)} = D_{\sqrt{c}}\;\rho_b \boxplus D_{\sqrt{1-c}}\;\gamma_0.
\end{equation}
\end{proposition}
This is a simple case of the following reformulation of the Main Theorem:

\begin{theorem} If $\mu$ is \emph{symmetric} probability measure on $\mathbb{R}$ with all moments, such that the even moments of the measure $\mu$ are following:
$$
m_{2n}(\mu) = \sum_{V \in \st{P}_2(2n)} t(V),
$$
and $t$ is normalized and \emph{strongly multiplicative}, then for $0 \leq b \leq 1$, the moments of the measure $D_{\sqrt{b}}(\mu) \boxplus D_{\sqrt{1-b}}\;\gamma_0 = \mu_b$, are of the form:
$$
m_{2n}(\mu_b) = \sum_{V \in \st{P}_2(2n)} b^{H(V)} \cdot t(V).
$$
and the \emph{free cumulants} of the measure $\mu$ and $\mu_b$ are following (see [B-Sp1] and [B-D-J], page 96):
$$
r_{2n}(\mu) = \hspace*{-3mm}\sum_{\scriptsize
\begin{array}{c}
V \in \st{P}_2(2n)\\
cc(V)=1\\
\end{array}
} \hspace*{-5mm}t(V),
$$
$$
r_{2n}(\mu_b)=b^n \hspace*{-5mm}\sum_{\scriptsize
\begin{array}{c}
cc(V)=1\\
V \in \st{P}_2(2n)\\
\end{array}
} \hspace*{-5mm}t(V) = b^n\; r_{2n}(\mu),\quad \mbox{for \underline{$n>1$}, $0 \leq b \leq 1$}.
$$
and $r_2(\mu_b)=r_2(\mu)$.
\end{theorem}

\begin{remark}
If a measure $\mu$ is free infinitely divisible (i.e. $\mu \in ID(\boxplus)$), then for $0 \leq b \leq 1$
$$
\mu_b = D_{\sqrt{b}}(\mu) \boxplus D_{\sqrt{1-b}}(\gamma_0) \mbox{ is also infinitely divisible},
$$
since Wigner semicircle law $\gamma_0 \in ID(\boxplus)$.

And vice verse,

\noindent if $\mu_b \in ID(\boxplus)$, for all $0 < b < 1$, then $\mu \in ID(\boxplus)$.

\end{remark}

\noindent {\bf Problem 2.} If we take $t(V) \equiv 1$,i.e. $\mu = \gamma_1$, $0 \leq b \leq 1$, then
$$
\sum_{V \in \st{P}_2(2n)} b^{H(V)} = m_{2n}(\omega_b)
$$
is the moment sequence of the probability measure $\omega_b = D_{\sqrt{b}}(\gamma_1) \boxplus D_{\sqrt{1-b}}(\gamma_0)$.

Is for $b>1$, the sequence
$$
m_{2n}(\mu) = \sum_{V \in \st{P}_2(2n)} b^{H(V)}
$$
the moment sequence of some symmetric probability measure?

See the paper [Boz1] on similar results for $q^{cr(V)}$, for $q>1$.

\section{Positive positive definite functions and ,,norm'' on permutation group}

In the paper [B-Sp2] we proved (Theorem~1) that if $t$ is positive definite function on $\st{P}_2(\infty)=\bigcup_{n=0}^\infty \st{P}_2(2n)$, than for all natural $n$, the restriction of $t$ to the permutation group $S(n)$ is also positive definite (in the usual sense), where the embedding $j:S(n) \rightarrow \st{P}_2(2n)$ is done later.

We recall that $t:\st{P}_2(\infty) \rightarrow \mathbb{C}$, is \emph{positive definite function} on $\st{P}_2(\infty)$, if there exists a~\emph{generalized Gaussian process (field)} $\{G_t(f)$, $f \in \st{H}\}$, $\st{H}$ -- real Hilbert space, such that
$$
\varepsilon [ G_t(f_1) G_t(f_2) \ldots G_t(f_k)]=\left\{
\begin{array}{ll}
0, & \mbox{if $k$ -- odd},\\
\displaystyle\sum_{V \in \st{P}_2(2n)} t(V) \displaystyle\prod_{(i,j) \in V} <f_i|j_j>, & k=2n\\
\end{array}
\right.
$$
for some state $\varepsilon$ on the $\ast$-algebra generated by the $G_t(f)$, $f \in \st{H}$ (see [G-M],[B-Sp2],[B-G],[Boz1],[B-Y], [B-W] for more examples of positive definite functions).

The our main theorem implies that the function $t_b(V)=b^{H(V)}$, for $0 \leq b \leq 1$, is positive definite on $\st{P}_2(\infty)$, which gives another proof of Buchholtz theorem \cite{Buch}.\\

\noindent Now we define the {\bf embedding} map $j:S(n) \rightarrow \st{P}_2(2n)$ formulated as follows:\\
for $\sigma \in S(n),\; j(\sigma)=\{(k,2n+1 - \sigma(k)): k=1,2,3,\ldots,n\} \in \st{P}_2(2n).$
On the picture:
\begin{center}
\includegraphics[width=7cm]{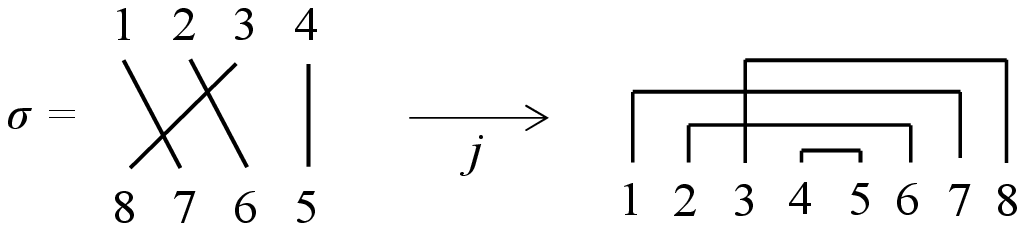}
\end{center}
From that figure we can see that the number of \emph{singletons} $h_n(\sigma) \stackrel{def}= h(j(\sigma))$, $\sigma \in S(n)$, is exactly the number of \emph{fixed points} of the permutation $\sigma$, which are \emph{isolated}.

The following Theorem is true:

\begin{theorem}\label{th5}
The function $h_{n+1}$ on $S(n+1)$ is of the form
$$
h_{n+1}=\sum_{k=1}^{n+1} s_{k-1} \cdot \tilde{s}_{(n+1)-k},
$$
and it is positive definite on the permutation group $S(n+1)$.
\end{theorem}

\noindent{\bf Proof of the Theorem~\ref{th5}.} Let $s_{k-1}=\chi_{s_{k-1}}$ is the characteristic function of the permutation group $S(k-1)$ on  $\{1,2,\ldots,k-1\}$, and
$\tilde{s}_{(n+1)-k}$ is the characteristic function of the symmetric group on the letters $\{k+1,k+2,\ldots,$ $n+1\}$,
\begin{center}
\includegraphics[width=4cm]{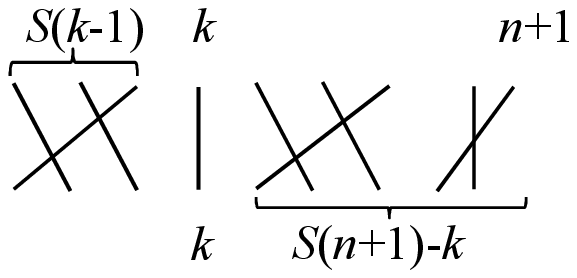}
\end{center}
here $S(k-1)$ is the group generated by inversions $\{\pi_1,\pi_2,\ldots,\pi_{k-2}\}$, and the group
$\tilde{S}(n+1-k)$ is generated by inversions $\{\pi_{k+1},\pi_{k+1},\ldots,\pi_{n}\}$
(here $\pi_j=(j,j+1)$ is the inversion (transposition) of $(j,j+1)$. Hence $s_{k-1} \cdot \tilde{s}_{(n+1)-k}$ is the cha\-ra\-cte\-ris\-tic function of the
group generated by $\{\pi_{1},\ldots,\pi_{k-2},\pi_k,\pi_{k+1},\ldots,\pi_{n}\}$, so it is positive definite on the group $S(n+1)$.

Let us define the function
$$h_{n+1}^{(k)}(\sigma)=\left\{
\begin{array}{ll}
1, & \mbox{if the singleton $(k,k)$ appears in the permutation $\sigma$},\\
0, & \mbox{otherwise}.\\
\end{array}
\right.
$$
\begin{center}
\includegraphics[width=28mm]{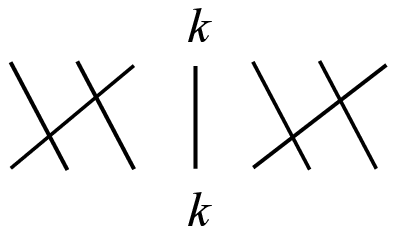}
\end{center}
For the above picture on can see that: $h_{n+1}^{(k)}=s_{k-1} \cdot \tilde{s}_{(n+1)-k}$.

Since our function
$$
h_{n+1}=\sum_{k=1}^{n+1} h_{n+1}^{(k)}, \mbox{ so we get}
$$
$$
h_{n+1}=\sum_{k=1}^n s_{k-1} \cdot \tilde{s}_{(n+1)-k}
$$
and it is positive definite as finite sum of positive definite functions.\hspace{5mm}\blacksquare

%Now we can state:

\begin{corollary}
For each $b \geq 1$ and $\sigma \in S(n)$
\begin{enumerate}
\item[($\alpha$)] The function $S(n) \ni \sigma \rightarrow b^{h_n(\sigma)}$ is positive definite on $S(n)$.

\item[($\beta$)] The function $H_n(\sigma)=n-h_n(\sigma)$ is conditionally negative definite on $S(n)$
(i.e. $\exp(-x\; H_n(\sigma))$ is positive definite on $S(n)$ for all positive $x>0$).

\item[($\gamma$)] The function $H_n(\sigma)$ is well defined on $S(\infty)=\bigcup S(n)$, $S(n) \subset S(n+1)$
(natural embedding) and $H_n=H_{n+1}|{S(n)}$, so we can define $H:S(\infty) \rightarrow \mathbb{R}$, as
$$
H(\sigma)=H_n(\sigma)=n-h_n(\sigma), \mbox{ for $\sigma \in S(n)$.}
$$
\end{enumerate}
\end{corollary}

\noindent {\bf Proof.} The case $(\gamma)$ can be easily checked, since the function\\
$h_{n+1}^{(k)}$ is the characteristic function of the group generated by
$$
\{\pi_1, \pi_2,\ldots,\pi_{k-2},\pi_{k+1},
\pi_{k+2},\ldots,\pi_{n}\},
$$
so by the restriction of
$h_{n+1}^{(k)}$ to the subgroup $\{\pi_1,\ldots,\pi_{k-2},\pi_{k+1},\ldots,\pi_{n-1}\}$,
we get
$h_{k+1}^{(k)}=h_n^{(k)}+1$, so
$$
H_n=H_{n+1}|{S(n)}.
$$

Both cases $(\alpha)$ and $(\beta)$ follow from the theorems of I. Schur and I. Schoenberg (see [Boz0] or Berg, Christensen, Ressel [B-Ch-R]).

Now we can state:

\begin{theorem}
The function $H$ is a ,,norm'' on $S(n)$ and also on $S(\infty)$, i.e.
\begin{enumerate}
\item[(i)] $H(e)=0$
\item[(ii)] $H(\sigma)=H(\sigma^{-1})$,\; $\sigma \in S(\infty)$
\item[(iiii)] $H(\sigma^{-1} \tau) \leq H(\sigma)+H(\tau)$, \; $\sigma,\tau \in S(\infty)$
\item[(iv)] If we define a function $d(\sigma,\tau)=H(\sigma^{-1} \tau)$, \; $\sigma,\tau \in S(\infty)$, then
$d$ is a left-invariant \underline{metric} on the group $S(\infty)$.

\end{enumerate}
\end{theorem}

\noindent {\bf Proof.} Let us see that for for $\sigma \in S_n$
$$
H(\sigma)=n-\sum_{k=0}^{n-1} s_k(\sigma) \cdot \tilde{s}_{n-k}(\sigma)=\sum_{k=0}^{n-1} (1-\chi_{s_k \times \tilde{s}_{n-k-1}}).
$$
Let us denote $\Delta_k=1-\chi_{s_k \times \tilde{s}_{n-k-1}}$, $k=0,1,2,\ldots,n-1$, then $\Delta_k$ is conditionally negative definite and
$$
\Delta_k(\sigma)=\left\{
\begin{array}{ll}
0, & \sigma \in S_k \times \tilde{S}_{n-1-k}\\
1, & \mbox{otherwise,}\\
\end{array}
\right.
$$
therefore $\sqrt{\Delta_k}=\Delta_k$.

By the well known property of the conditionally negative definite function (see [B-Ch-R], [Boz0]), we have
$$
\Delta_k(\sigma \tau)=\sqrt{\Delta_k(\sigma \tau)} \leq \sqrt{\Delta_k(\sigma)} + \sqrt{\Delta_k(\tau)} = \Delta_k(\sigma) + \Delta_k(\tau).
$$
Therefore our function
$$
H=\sum_{k=0}^{n-1} \Delta_k
$$
is also subadditive.

Other facts follow at once from the definition of the function H.\;\blacksquare\\

\section{Some questions and problems}

\noindent {\bf Problem 3.} If the assumption that the function $t:\st{P}_2(\infty) \rightarrow \mathbb{R}$ is \underline{strongly} multiplicative
is necessary in our Main Theorem and Lemmas?\\

\noindent {\bf Problem 4.} Let $\Gamma_b(\st{H})$ is the von Neumann algebra generated by our $b$-Gaussian process
$$
Y_b(f)=\sqrt{1-b} G_0(f) + \sqrt{b} G_1(f), \quad f \in \st{H}_{\mathbb{R}},
$$
where $G_0 (f)$ is the free Gaussian process and $G_1 (f)$ is the classical Gaussian process and $G_0(f)$ and $G_1(f)$ are free independent.
As we can see in that von Neumann algebra $\Gamma_b (\st{H})$, the vacuum is the trace, it is faithful and normal state.\\

\noindent {\bf Problem 5.} Natural question is, if that von Neumann algebras
$\Gamma_b(\st{H})$ is a factor, $\Gamma_b(\st{H}) = \{ Y_b(f):f \in \st{H}_{\mathbb{R}}\}''$ (bicommutant), for $\rm{dim}\; \st{H} \geq 2$.

\noindent {\bf Problem 6.}  If for $0 < b < 1$ that algebras
$\Gamma_b(\st{H})$ are isomorphic as von Neumann algebras?

%One can prove that algebras $\Gamma_b(\st{H})$ are \underline{not injective} as von Neumann algebras, if $0 \leq b < 1$, and $\rm{dim}\; \st{H} \geq %2$.

Some others facts about our algebra $\Gamma_b(\st{H})$ will be presented in the second part of our paper.

\section*{Acknowledgement}

The work was partially supported by the MAESTRO grant DEC-2011/02/A/\\ST1/00119 and OPUS grant DEC-2012/05/B/ST1/00626 of National Centre of Science (M. Bo\.{z}ejko) and the OPUS grant DEC-2012/05/B/ST7/00102 of National Centre of Science (W. Bo\.{z}ejko).

\end{document}